\documentclass[12pt]{article}
\usepackage{amssymb}
\usepackage{latexsym}

\oddsidemargin=\evensidemargin
\begin{document}
\newfont{\blb}{msbm10 scaled\magstep1} 

\newtheorem{defi}{Definition}
\newtheorem{theo}{Theorem}[section]
\newtheorem{theo1}{Theorem}[section]
\newtheorem{prop}[theo]{Proposition}
\newtheorem{lemm}[theo]{Lemma}
\newtheorem{coro}[theo]{Corollary}
\date{}
\author{G. Traustason,
University of Bath, UK \\ 
J. Williams, University of Bristol, UK} 
%
\title{Powerfully nilpotent groups  of rank 2 or small order}
\maketitle
\begin{abstract}
In this paper we continue the study of powerfully nilpotent groups started in [3] and continued in [4,5,6]. These are powerful  $p$-groups possessing a central series of a special
kind. To each such group one can attach a powerful nilpotency class that leads naturally to the notion of a powerful coclass and classification in terms of an ancestry
tree. In this paper we will give a full classification of powerfully nilpotent groups of rank $2$. The classification will then be used to arrive at a precise formula for the number of powerfully nilpotent groups of rank $2$ and order $p^{n}$. We will also give a detailed analysis of
the ancestry tree for these groups. The second part of the paper is then devoted to a full classification of powerfully nilpotent groups
of order up to $p^{6}$. 
\end{abstract}
\section{Introduction}
In this paper we continue the study of powerfully nilpotent $p$-groups started in [3] and continued in [4,5,6]. 
Powerful $p$-groups were introduced by Lubotzky and Mann in [2].  The class of powerfully nilpotent groups is
a special subclass of these, containing groups that possess a central series of a special kind. 
We start by recalling
the basic terms. Let $G$ be a finite $p$-group where $p$ is a prime. \\ \\
{\bf Definition}. Let $H\leq K\leq G$. An ascending chain of subgroups 
  $$H=H_{0}\leq H_{1}\leq \cdots \leq H_{n}=K$$
is {\it powerfully central in $G$} if $[H_{i},G]\leq H_{i-1}^{p}$ 
for $i=1,\ldots ,n$. Here $n$ is called the {\it length} of the chain. \\ \\
{\bf Definition}. A powerful $p$-group 
 $G$ is {\it powerfully nilpotent} if it has an ascending chain of subgroups of the form
    $$\{1\}=H_{0}\leq H_{1}\leq \cdots \leq H_{n}=G$$ 
that is powerfully central in $G$. \\ \\
{\bf Definition}. If $G$ is powerfully nilpotent then the {\it powerful nilpotence class} of $G$ is the shortest length that a powerfully central chain of $G$ can have. \\ \\
{\bf Definition}. We say that a finite $p$-group $G$ is {\it strongly powerful} if $[G,G]\leq G^{p^{2}}$. \\ \\
In [6] we showed that a strongly powerful $p$-group is powerfully nilpotent of powerful class at most $e-1$ where
$p^{e}$ is the exponent of the group. \\ \\
{\bf The upper powerfully central series}. This is defined recursively as 
follows: $\hat{Z}_{0}(G)=\{1\}$ and for $n\geq 1$, 
      $$\hat{Z}_{n}(G)=\{a\in G:\, [a,x]\in \hat{Z}_{n-1}(G)^{p}\mbox{ for
all }x\in G\}.$$
Notice in particular that $\hat{Z}_{1}(G)=Z(G)$. \\ \\
{\bf Definition}. Let $G$ be a powerfully nilpotent $p$-group of powerful
class $c$ and order $p^{n}$. We define the {\it powerful coclass} of $G$ to be the number $n-c$. \\ \\
A natural approach is to develop something that corresponds to a coclass theory for finite $p$-groups where coclass is replaced by powerful coclass and in [3] we proved that there are indeed, for any fixed prime $p$, finitely many powerfully nilpotent $p$-groups of any given powerful coclass.  More precisely, if $G$ is
a powerfully nilpotent $p$-group of rank $r$ and exponent $p^{e}$ then we showed that  
$r\leq n-c+1$ and $e\leq n-c+1$. This together with the fact that $n\leq re$ in a powerful group gives the result. \\ \\
{\bf The ancestry tree}. Let $p$ be a fixed prime. The {\it vertices} of the 
ancestry tree are all the powerfully nilpotent $p$-groups (one for each isomorphism class). Two vertices $G$ and $H$ are joined by a directed edge from $H$ to
$G$ if and only if $H\cong G/Z(G)^{p}$ and $G$ is not abelian. Notice that this implies that $Z(G)^{p}\not =\{1\}$ and thus the powerful class of $G$ is one more than that of $H$. 
We then also say that 
$G$ is a {\it direct ancestor} of $H$ or that $H$ is a {\it direct descendant}
of $G$, and we write $H\rightarrow G$.   \\ \\ 
Let $c(H),\,c(G)$ be the powerfully nilpotent classes of $H$ and $G$ and let $d(H),\,d(G)$ be the powerful coclasses, that is $d(H)=n(H)-c(H)$ and $d(G)=n(G)-c(G)$ where the orders of $H$ and $G$ are $p^{n(H)}$ and
$p^{n(G)}$. Notice that $c(H)=c(G)+1$ and that $d(H)\geq d(G)$ with equality if and only if $|Z(G)^{p}|=p$.  \\ \\
%
%
%
%
%
%
%
%
%
%
%
%
Let $G$ be a powerfully nilpotent $p$-group of rank $r$ and exponent $p^{e}$. 
As $G$ is powerful we can choose our generators $a_{1},\ldots ,a_{r}$ such that $G=\langle a_{1}\rangle\cdots \langle a_{r}\rangle$
and such that $|G|=o(a_{1})\cdots o(a_{r})$. In [3, Theorem 2.7] we furthermore showed that the generators can be chosen such that the following chain is powerfully central
$$ \begin{array}{c}
    G\geq \langle a_{1}^{p},a_{2},\ldots ,a_{r}\rangle\geq \langle a_{1}^{p},a_{2}^{p},a_{3},\ldots ,a_{r}\rangle \geq\cdots\geq \langle a_{1}^{p},a_{2}^{p},\ldots ,a_{r}^{p}\rangle = \\
   G^{p}\geq \langle a_{1}^{p^{2}},a_{p},\ldots ,a_{r}^{p}\rangle\geq \langle a_{1}^{p^{2}},a_{2}^{p^{2}},a_{3}^{p},\ldots ,a_{r}^{p}\rangle \geq\cdots\geq \langle a_{1}^{p^{2}},a_{2}^{p^{2}},\ldots ,a_{r}^{p^{2}}\rangle = \\
                                            \vdots \\
    G^{p^{e-1}}\geq \langle a_{1}^{p^{e}},a_{p},\ldots ,a_{r}^{p}\rangle\geq \langle a_{1}^{p^{e}},a_{2}^{p^{e}},a_{3}^{p},\ldots ,a_{r}^{p}\rangle \geq\cdots\geq \langle a_{1}^{p^{e}},a_{2}^{p^{e}},\ldots ,a_{r}^{p^{e}}\rangle = 1.
\end{array}$$
One can furthermore pick the generators such that those of order $p$ come before the others. If there are $n_{k}$ generators of
order $p^{k}$, then we say that the group is of type 
                 $$(\underbrace{1,\ldots ,1}_{n_{1}},\cdots ,\underbrace{e,\ldots ,e}_{n_{e}}).$$
Sometimes we use the shorter notation $(1_{n_{1}},\ldots ,e_{n_{e}})$ for this. \\ \\
This is something we will make much use of in the classification below. \\ \\
The paper is organised as follows. In Section 2 we will give a complete classification of all powerfully nilpotent groups of rank $2$ and then use that classification to obtain a precise closed formula for the number of powerfully nilpotent groups of rank $2$ that
are of order $p^{n}$.  In Section 3 we will then give a detailed analysis of the ancestry tree of powerfully nilpotent $p$-groups of rank $2$. In the final section we will give a complete classification of all powerfully nilpotent groups of order up to $p^{6}$.
\section{Powerfully nilpotent groups of rank $2$: Enumeration}
In this section we give, for any given prime $p$, a full classification of the $2$-generator powerfully nilpotent $p$-groups. Using the classification we obtain a closed formula for the number of powerfully nilpotent $2$-generator groups of order $p^{n}$. \\ \\
It is well known that all powerful $2$-generator $p$-groups are metacyclic. In fact if $G=\langle a,b\rangle$ is a $2$-generator powerful $p$-group, then we know that $[G,G]=\langle [a,b]\rangle^{G}$ is powerfully embedded in $G$ and thus $[G,G]=\langle [a,b]\Phi([G,G])\rangle =\langle [a,b]\rangle$. As $G$ is powerful
we have $[a,b]=c^{p^{r}}$ for some $r\geq 1$ and $c\in G\setminus \Phi(G)$. It follows that we can without loss of generality assume that $[a,b]\in \langle a^{p^{r}}\rangle$ and thus $\langle a\rangle$ is normal in $G$ from which it follows that $G$ is metacyclic. \\ \\
There is a classification of metacyclic $p$-groups [1] and in principle one could work out from this a classification of all powerfully nilpotent $p$-groups of rank $2$. We prefer though to give
a short direct proof that reflects better the structure of powerfully nilpotent groups and thus provides better insight. We start with a simple observation. 
\begin{lemm} Any $2$-generator powerfully nilpotent $p$-group is strongly powerful.
\end{lemm}
{\bf Proof}\ \ Let $G$ be a $2$-generator group that is powerfully nilpotent. If $p=2$, then it follows straight from the definition that $G$ is strongly powerful. We
can thus assume that $p$ is odd. By [3, Theorem 2.7] we can pick generators $a,b$ for $G$ such that the chain
       $$G=\langle a,b\rangle \geq \langle a^{p},b\rangle \geq \langle a^{p},b^{p}\rangle=G^{p}$$
is powerfully central. In particular it follows from this that, for $H=\langle a^{p},b\rangle$, we have $[G,G]=[H,G]\leq (G^{p})^{p}=G^{p^{2}}$. $\Box$ \\ \\
Let $G$ be a non-abelian powerfully nilpotent group of rank $2$ and let $r$ be the largest positive integer such that $[G,G]\leq G^{p^{r}}$. By Lemma 2.1 we know that $r\geq2$. Then there exists $n>r$
such that $|[G,G]|=p^{n-r}$. As $[G,G]$ is a cyclic subgroup of
$G^{p^{r}}$ and $G$ is powerful, $[G,G]$ must have a generator of
the form $a^{p^{r}}$ for some $a\in G$. By the choice of $r$ we must have that $a\not\in G^{p}$. Notice that $o(a)=p^{n}$ and that $\langle a\rangle$
is normal in $G$ as $G'\leq\langle a\rangle$. Let $m$ be a positive
integer such that $p^{m}=\left|\frac{G}{\langle a\rangle}\right|=\frac{|G|}{p^{n}}.$
It follows that $|G|=p^{n+m}$. For any $b\in G\backslash\langle a\rangle G^{p}$  we have $b^{p^{m}}\in\langle a\rangle$,
say $b^{p^{m}}=a^{p^{l}},$ with $1\leq l\leq n$. Among all choices
of $a$ and $b$ as above, make a choice so that $l$ is maximal.
With a slight adjustment of the choice of generators (left to the reader) one may then assume we have a presentation of the form: $a^{p^{n}}=1$,
$b^{p^{m}}=a^{p^{l}}$ and $[a,b]=a^{p^{r}}$. From this one can read that  $Z(G)=\langle a^{p^{n-r}},b^{p^{n-r}}\rangle$. Hence
%
 we must have $m,l\geq n-r$, as $b^{p^{m}},a^{p^{l}}\in Z(G)$. \\ \\
Notice also that either $l=n$ in which case the group is a semidirect
product of a cyclic group of order $p^{n}$ by a cyclic group of order $p^{m}$, or else we must have $m>l$. To see this, suppose for contradiction
that $m\leq l<n$. Modulo $\langle a^{p^{l+1}}\rangle$, we then have $(ba^{-p^{l-m}})^{p^{m}}=b^{p^{m}}a^{-p^{l}}=1$
and thus $(ba^{-p^{l-m}})^{p^{m}}\in\langle a^{p^{l+1}}\rangle$ contradicting
the maximality of $l$. Finally, for the non-split groups we must
have $r<l$. To see this, suppose for contradiction that $l\leq r$. Then $a^{p^{r}}=(a^{p^{l}})^{p^{r-l}}=(b^{p^{m}})^{p^{r-l}}=b^{p^{m+r-l}}$ from which it would then follow that $[a,b]=b^{p^{m+r-l}}$. However, as $m+r-l>r$ this would imply that $[G,G]\leq G^{p^{r+1}}$ contradicting the choice of $r$. From this analysis we arrive at two types of presentations: \\ \\
\underline{I. Semidirect products:} $G(n,m,r):=\langle a,b|a^{p^{n}}=1,b^{p^{m}}=1,[a,b]=a^{p^{r}}\rangle$
where the  parameters satisfy the  inequalities  \\
$$\begin{array}{l}
n-r \leq m \\
2\leq  r\leq n-1.
\end{array}$$
\underline{II. Non-semidirect products:} $G(n,m,l,r):=\langle a,b|a^{p^{n}}=1,b^{p^{m}}=a^{p^{l}},[a,b]=a^{p^{r}}\rangle$
with $4$ parameters satisfying the inequalities 
%
%
%
$$\begin{array}{l}
2\leq  r< l\,\leq  n-1\\
\mbox{}n-r \leq  l < m.
\end{array}$$
{\bf Remark}. These presentations are consistent and do define groups of order $p^{n+m}$.
For groups of type I this is clear since they are semidirect products and $o(b)=p^{m}$ is consistent with $a^{b}=a^{1+p^{r}}$ as $m+r\geq n$. 
For groups of type II consider first the group $H$ with presentation $\langle\bar{a},\bar{b}|\bar{a}^{p^{n}}=1,\bar{b}^{p^{m+n-l}}=1,[\bar{a},\bar{b}]=\bar{a}^{p^{r}}\rangle$, and note that the presentation for $H$ is of type I. Consider the subgroup $\langle\bar{b}^{p^{m}}\bar{a}^{-p^{l}}\rangle$.
Note that as $m,l\geq n-r$ we have that $\bar{b}^{p^{m}}\bar{a}^{-p^{l}}\in Z(H)$
and so the subgroup is normal. Then $G$ can be realised as the quotient
$H/\langle\bar{b}^{p^{m}}\bar{a}^{-p^{l}}\rangle$.\\ \\
{\bf Remark}. Notice that for $G=G(n,m,r)$ or $G=G(n,m,l,r)$, $r$ is a structural invariant determined from $r$ being the largest positive integer such that $[G,G]\leq G^{p^{r}}$ (for the latter we need $r<l$). This implies that $n$ and $m$ are also structure invariants as $|[G,G]|=p^{n-r}$ and $|G|=p^{n+m}$. Next notice that if $G=G(n,m,r)$ then $|G^{p^{m}}|=1$ if $m\geq n$ and $|G^{p^{m}}|=o(a^{p^{m}})=p^{n-m}$ if $m<n$. On the other hand if $G=G(n,m,l,r)$ then $|G^{p^{m}}|=|\langle a^{p^{m}},b^{p^{m}}\rangle|=|\langle a^{p^{l}}\rangle|=p^{n-l}$ as $l<m$. Notice also that $n-l>n-m$. This shows that the groups of type II are distinct from 
the groups of type I and that for groups of type II, $l$ is a structure invariant determined from $|G^{p^{m}}|=p^{n-l}$. It follows from this analysis that there is no duplication in the list of the groups of type I and II and thus we have a complete classification of the powerfully nilpotent $2$-generator groups of rank $2$. \\ \\
We next use this classification to enumerate the powerfully nilpotent $p$ groups of rank $2$ and order $p^{x}$. With respect to the notation above we have $x=n+m$. Let $y=n-r$. With respect to the variables $x,y,n$ the system of inequalites for groups of type I is equivalent to: 
$$\begin{array}{ll} 
 & 1\leq y\leq n-2 \\
 & y+(n-2)\leq x-2.
\end{array}$$
This is equivalent to $1\leq y\leq \lfloor\frac{x-2}{2}\rfloor$ (where $\lfloor z\rfloor$ is the floor of $z$) and $y+2\leq n\leq x-y$. We thus need to find the number of all integer pairs $(n,y)$ satisfying this system of inequalities. Notice that for a solution to exist we need $1\leq \lfloor\frac{x-2}{2}\rfloor$, that is $x\geq 4$. For a given $x\geq 4$ the number of solutions 
is
\begin{equation}
\small
 \sum_{y=1}^{\lfloor \small \frac{x-2}{2}\rfloor }(x-1-2y)= \left[ \begin{array}{l} (\frac{x-2}{2})^{2} \mbox{\ \ \ if }x\mbox{ is even,}\\
                                     \frac{(x-3)(x-1)}{4} \mbox{\ if }x \mbox{ is odd}. \end{array}\right.
\end{equation}
\normalsize
We then move on to groups of type II. Working again with the parameters $x,y,n$ as well as $l$ we get the following equivalent system
of equations
$$\begin{array}{l}
l+1\leq n\leq x-1-l \\
n+1-l\leq y\leq l \\
2\leq n-y
\end{array}$$
Notice that it follows from these inequalities that $2\leq n-y\leq l-1$ and $l+1\leq x-1-l$ from which we get that $3\leq l\leq \lfloor\frac{x-2}{2}\rfloor$. For a solution 
to exist we thus need $x\geq 8$.  Notice also that if $n>l+1$ then from $y\leq l$ it follows that $n-y\geq l+2-l=2$ and thus the last inequality will be a consequence of the others. We thus deal first
separately with the case $n=l+1$. \\ \\
{\bf Case (a).} Suppose $n=l+1$. We then get the inequalites $2\leq y\leq l$ and $y\leq l-1$. Thus for a given $l$ such that $3\leq l\leq \lfloor\frac{x-2}{2}\rfloor$
we get \underline{$l-2$} solutions $(n,y)=(l+1,y)$. \\ \\
{\bf Case (b)}. Suppose $l+2\leq n\leq x-1-l$. We need to count the solutions to  
$$\begin{array}{l}
l+2\leq n\leq x-1-l \\
n+1-l\leq y\leq l .\\
\end{array}$$
Notice that we need $x\geq 9$ for a solution to exist (for $x=8$ we have $n=l+1$). For the second set of inequalities to hold we need $n\leq 2l-1$. When $x< 3l$ then $n\leq x-1-l\leq 2l-1$ and thus
for all $l$ such that $\lfloor\frac{x}{3}\rfloor<l\leq \lfloor \frac{x-2}{2}\rfloor$ there are $2l-n$ solutions for each $n$ where $l+2\leq n\leq x-1-l$. 
 On the other hand
when $3\leq l\leq \lfloor\frac{x}{3}\rfloor$, then $2l-1\leq x-1-l$ and solutions only exists for $n$ when $l+2\leq n\leq 2l-1$. \\ \\
Adding up, the number of groups of type II is
$$\sum_{l=3}^{\lfloor\frac{x}{3}\rfloor}\left[\left(\sum_{n=l+2}^{2l-1}(2l-n)\right)+(l-2)\right] + \sum_{l=\lfloor\frac{x}{3}\rfloor+1}^{\lfloor\frac{x-2}{2}\rfloor}\left[\left(\sum_{n=l+2}^{x-1-l}(2l-n)\right)+(l-2)\right]$$
After tedious but straightforward calculations we see that this equals 
 $$
\begin{array}{r}\left\lfloor \frac{x-6}{2}\right\rfloor \frac{(x-6)(7-x)}{2}+\left\lfloor \frac{x-6}{2}\right\rfloor \left\lfloor \frac{x-4}{2}\right\rfloor \frac{(3x-18)}{2} -  \frac{4\lfloor\frac{x-6}{2}\rfloor\lfloor\frac{x-4}{2}\rfloor(2\lfloor\frac{x-6}{2}\rfloor+1)}{6}+\left\lfloor \frac{x-6}{3}\right\rfloor \frac{(x-6)
(x-7)}{2} \\
  +  \left\lfloor \frac{x-6}{3}\right\rfloor \left\lfloor \frac{x-3}{3}\right\rfloor \left(\frac{-6x+39}{4}\right)+\frac{9}{2}\frac{\lfloor\frac{x-6}{3}\rfloor\lfloor\frac{x-3}{3}\rfloor(2\lfloor\frac{x-6}{3}\rfloor+1)}{6}.
\end{array}$$


\mbox{}\\ 
Adding this to the number of groups of type 1 as well as those groups that are abelian of rank $2$ (whose number is $\lfloor\frac{x}{2}\rfloor$)
gives us the following result.
\begin{prop}
For $x\geq 4$, the number of rank $2$ groups of order $p^{x}$ which are powerfully nilpotent is 
\small
\begin{eqnarray*}
\frac{x^{3}+12x^{2}-60x+216}{72} & \,\mbox{if } x\equiv 0(\mbox{mod }6) \\
\frac{x^{3}+12x^{2}-69x+200}{72} & \,\mbox{if } x\equiv 1(\mbox{mod } 6)\\
\frac{x^{3}+12x^{2}-60x+208}{72} & \,\mbox{if } x\equiv 2(\mbox{mod }6)\\
\frac{x^{3}+12x^{2}-69x+216}{72} & \,\mbox{if } x\equiv 3(\mbox{mod }6)\\
\frac{x^{3}+12x^{2}-60x+200}{72} & \,\mbox{if }x\equiv 4(\mbox{mod }6)\\
\frac{x^{3}+12x^{2}-69x+208}{72} & \,\mbox{if }x\equiv 5(\mbox{mod }6)
\end{eqnarray*}
\normalsize
\end{prop}

\section{The ancestry tree for groups of rank 2}

As a second application of the complete classification of powerfully nilpotent groups of rank $2$, we will in this section give a detailed analysis of the ancestry tree. 
We start by determining the direct descendants of the powerfully nilpotent groups of rank $2$.  We consider two separate cases. \\ \\
{\bf Direct descendants of the split groups}. Let $G=\langle a,b:\,a^{p^{n}}=1,\ b^{p^{m}}=1\, [a,b]=a^{p^{r}}\rangle$ be a $2$-generator powerfully nilpotent group that is a semidirect product of a cyclic group of order $p^{n}$ by cyclic group of order $p^{m}$. According to the classification the three parameters $m,n,r$ satisfy 
                $$m\geq n-r\mbox{\ and }2\leq r\leq n-1.$$
We have seen above that $Z(G)^{p}=\langle a^{p^{n-r+1}}, b^{p^{n-r+1}}\rangle$. There are now two subcases to be considered depending on whether
$m>n-r$ or $m=n-r$. Now the direct descendant of $G$ is $D(G)=G/Z(G)^{p}$. Thus we have \\ \\
(a) $D(G)=\langle a,b:\, a^{p^{n-1+1}}=1,\ b^{p^{n-r+1}}=1,\ [a,b]=a^{p^{r}}\rangle$ if $m>n-r$. \\
(b) $D(G)=\langle a,b:\,  a^{p^{n-1+r}}=1,\ b^{p^{n-r}}=1,\ [a,b]=a^{p^{r}}\rangle$ if $m=n-r$. \\ \\ 
Notice that in both cases $D(G)$ is abelian if and only if $r\geq n-r+1$.  \\ \\
{\bf Direct descendants of the non-split groups}. Let $G=\langle a,b:\, a^{p^{n}}=1,\ b^{p^{m}}=a^{p^{l}},\ [a,b]=a^{p^{r}}\rangle$ be one of
the non-split groups. According to the classification we have
         $$2\leq r<l\leq n-1\mbox{ and }n-r\leq l<m.$$ 
As before we have two cases, this time depending on the relationship between $n-r$ and $l$. Notice that if $n-r=l$ then $Z(G)^{p}=\langle a^{p^{l+1}},b^{p^{l+1}}\rangle$. But as  $l+1\leq m$ we then have $b^{p^{m}}=a^{p^{l}}\in Z(G)^{p}.$ Thus we have \\ \\
(c) $D(G)=\langle a,b:\,a^{p^{n-r+1}}=1,\mbox{ }b^{p^{n-r+1}}=1,\mbox{ }[a,b]=a^{p^{r}}\rangle$ if $n-r<l$. \\
(d) $D(G)=\langle a,b:\,a^{p^{n-r}}=1,\mbox{ }b^{p^{n-r+1}}=1\mbox{, }[a,b]=a^{p^{r}}\rangle$ if $n-r=l$. \\ \\
Notice that in (c) we get an abelian group if and only if $r\geq n-r+1$. In (d) we never have an abelian group as for this to happen one would need $r\geq n-r$ but this cannot be the case as then $r\geq l$ contradicting the requirements for the parameters. \\ \\
{\bf Remark}. The analysis above shows that a direct descendant is always one of the split groups. In other words the non-split groups do not have
ancestors. \\ \\
Having determined the direct descendants, one can now use this to read what the direct ancestors are of any given group of rank $2$. We have already
seen that the non-split groups do not have ancestors and we thus only need to consider the split groups. \\ \\
{\bf Direct ancestors of abelian groups}. We see from the work above that for an abelian group to have an ancestor, it needs be of one of the following
forms: \\ \\
$A(\bar{n})=\langle a^{p^{\bar{n}}}=1,\ b^{p^{\bar{n}}}=1\rangle $, where $\bar{n}\geq 2$; \\
$B(\bar{n})=\langle a^{p^{\bar{n}+1}}=1,\ b^{p^{\bar{n}}}=1\rangle$, where $\bar{n}\geq 1$. \\ \\ 
For $G(n,m,r)$ or $G(n,m,l,r)$ to be a direct ancestor of $A(\bar{n})$  we need $\bar{n}=n-r$. The conditions on $\bar{n}$ then follow from the requirement
that $n-r\geq 1$. From the analysis of the direct descendants we can find all the direct ancestors of $A(\bar{n})$ and $B(\bar{n})$. \\ \\
\underline{The direct ancestors of $A(\bar{n})$}. The split ancestors are $G(\bar{n}+r-1,m,r)$ where $r\geq \bar{n}$ and $m\geq \bar{n}$. 
Here $r\geq \bar{n}$ since the the desecendant is abelian and $m\geq \bar{n}$ since this is the same as saying that
$m>(\bar{n}+r-1)-r$ which was one of the requirements for the parameters (see (a)). Notice that there are infinitely many such groups.  \\ \\
The non-split ancestors  are $G(\bar{n}+r-1,m,l,r)$ where $r+(\bar{n}-2)\geq l>r\geq \bar{n}$ and $m>l\geq \bar{n}$. Again the
condition $r\geq \bar{n}$ comes from the descendant. The rest comes from the requirements for the parameters. Notice that $l\geq \bar{n}$ is equivalent to $(\bar{n}+r-1)-r<l$ that comes from the conditition in (c).  \\ \\
\underline{The  direct ancestors of $B(\bar{n}$)}. In this case we can read from our analysis above that the ancestor must be a semidirect product $G(\bar{n}+r,\bar{n},r)$ with $r\geq \bar{n}+1$ (for the descendant to be abelian). \\ \\
{\bf The direct ancestors of the semidirect products}. In this case one can read that the semidirect products that have ancestors are: \\ \\
$C(\bar{n})=\langle a,b:\,a^{p^{\bar{n}}}=1,\ b^{p^{\bar{n}}}=1,\ [a,b]=a^{p^{r}}$ with $2\leq r\leq \bar{n}-1$ and $\bar{n}\geq 3$; \\
$D(\bar{n})=\langle a,b:\,a^{p^{\bar{n}+1}}=1,\ b^{p^{\bar{n}}}=1,\ [a,b]=a^{p^{r}}$ with $2\leq r\leq \bar{n}$ and $\bar{n}\geq 2$; \\
$E(\bar{n})=\langle a,b:\,a^{p^{\bar{n}}}=1,\ b^{p^{\bar{n}+1}}=1,\ [a,b]=a^{p^{r}}$ with $2\leq r\leq \bar{n}-1$ and $\bar{n}\geq 3$.  \\ \\
\underline{The direct ancestors of $C(\bar{n})$}. The semidirect ancestors are $G(\bar{n}+r-1, m, r)$ with $m\geq \bar{n}$ as this is the situation (a).  \\ \\
The non-split ancestors are $G(\bar{n}+1-1,m,l,r)$ with $m>l\geq \bar{n}$ and $l\leq r+\bar{n}-2$. For both these cases we get infinitely many direct ancestors. \\ \\
\underline{The direct ancestors of $D(\bar{n})$}. In this case we are in situation (b) and we know that all the direct ancestors will
be semidirect products and these are $G(\bar{n}+r,\bar{n},r)$. \\ \\
\underline{The direct ancestors of $E(\bar{n})$}. Here we are in situation (d) and we see that the direct ancestors are $G(\bar{n}+r,m,\bar{n},r$ with $m>\bar{n}$. Here again there are infinitely many direct ancestors. \\ \\
{\bf The infinite branches of the ancestry tree}. It follow from the analysis above that the only abelian groups that have grandparents
are $A(\bar{n})$ where $\bar{n}\geq 2$ and 
$G(\bar{n}, \bar{n},r)$ where $\bar{n}\geq 3$ and $2\leq r\leq \bar{n}-1$. \\ \\
In particular, if there is an infinite branch, it must start with one of the abelian groups $A(\bar{n})$ and then all the subsequent groups should be of the second type. We thus get a chain of ancestors of the following form:
    $$A(\bar{n})\rightarrow G(\bar{n}+r-1,\bar{n}+r-1,r)\rightarrow \cdots \rightarrow G(\bar{n}+i(r-1),\bar{n}+i(r-1),r)\rightarrow \cdots $$
where there is one such infinite branch for each pair $(\bar{n},r)$ where $\bar{n},r\geq 2$. 
\section{Classification of powerfully nilpotent $p$-groups of order less than or equal to $p^{6}$.}
In this final section we classify all powerfully nilpotent $p$-groups of order up to and including $p^{6}$. We will deal with cases when $p$ is odd and when $p=2$ separately. We start with the case when $p$ is odd.
\subsection{The case when $p$ is odd}
We first remark that there are no non-abelian $p$-groups of order less than $p^{4}$ that are powerfully nilpotent. To see this we argue  by contradiction and suppose we have a non-abelian powerfully nilpotent group $G$ of order $p^{3}$. Then we must have that $Z(G)^{p}\not=\{1\}$ and therefore that $Z(G)$ has order $p^{2}$. But in that case $G/Z(G)$ is cyclic that gives the contradiction that $G$ is abelian. We thus only need to deal with groups of order $p^{4}$, $p^{5}$ and $p^{6}$. For groups of rank $2$ we can read from Section 2 what the non-abelian groups are. We thus only need to do some work for the non-abelian groups of rank greater than or equal to $3$. \\ \\ 
The following general setting includes a number of groups that  occur, namely groups of type $(1_{t}, n)$ where $n$ is an integer greater than $1$.  Suppose $G=\langle a_{1},\ldots ,a_{t},b\rangle$ is a powerfully nilpotent group of this type where $a_{i}$ is of order $p$ and $b$ of order $p^{n}$. Notice that $G^{p}=\langle b^{p}\rangle$ is cyclic and it follows from [3, Corollary 3.3] that $G^{p}\leq Z(G)$.
In particular $G$ is nilpotent of class at most $2$. As a result $[G,G]^{p}=[G^{p},G]=1$. Next observe that $\Omega_{1}(G)=
\langle a_{1},\ldots ,a_{r},b^{p^{n-1}}\rangle$ where $\Omega_{1}(G)$ is the subgroup consisting of all elements of order dividing $p$. Thus $[G,G]\leq \Omega_{1}(G)\cap G^{p}=\langle a^{p^{n-1}}\rangle$. \\ \\
Consider the vector space $V=\Omega_{1}(G)G^{p}/G^{p}=\langle a_{1},\ldots ,a_{t}\rangle G^{p}/G^{p}$ over the field ${\mathbb F}$ of $p$ elements. The commutator operation induces naturally an alternating form on $V$ through
             $$(xG^{p},yG^{p})=\lambda\mbox{ if }[x,y]=a^{\lambda p^{n-1}}.$$
Without loss of generality we can suppose that our generators have been chosen such that we get the following orthogonal decomposition 
          $$V=\langle a_{1}G^{p}, a_{2}G^{p}\rangle\oplus \cdots \oplus \langle a_{2s-1}G^{p},a_{2s}G^{p}\rangle\oplus V^{\perp},$$
where $V^{\perp}=\langle b_{2s+1}G^{P},\ldots ,b_{t}G^{p}\rangle$ and $(a_{2i-1}G^{p},a_{2i}G^{p})= 1$ for $i=1,\ldots ,s$. There are now two cases to consider, depending on whether or not $Z(G)\leq \Omega_{n-1}(G)$, where $\omega_{n-1}(G)$ is the subgroup of elements of order dividing $p^{n-1}$. \\ \\
Suppose first  that $Z(G)\not\leq \Omega_{n-1}(G)$. This means that $Z(G)$ contains some element $b^{r}u$ with $u\in \langle
a_{1},\ldots ,a_{t}\rangle$ and $0<r<p$. One sees readily that one can choose the set of generators to include $b^{r}u$. Thus without loss of generality we can assume that $b\in Z(G)$. We thus get a powerfully nilpotent group $A(n,t,s)$ with powerfully nilpotent presentation that has generators $a_{1},\ldots ,a_{t},b$ and relations
$$\begin{array}{l}
                    a_{1}^{p}=\cdots =a_{t}^{p}=b^{p^{n}}=1, \\
                    \mbox{} [a_{2i-1},a_{2i}]=a^{p^{n-1}}\mbox{ for }i=1,\ldots ,s, \\
           \mbox{}        [a_{i},a_{j}]=1\mbox{ otherwise for }1\leq i<j\leq t, \\
      \mbox{}            [a_{i},b]=1\mbox{ for }1\leq i\leq t.
\end{array}$$
Notice that for a fixed $n\geq 2$ and $t\geq 2$ we get $\lfloor t/2\rfloor $ such non-abelian roups of order $p^{n+t}$, namely $A(n,t,s)$
for $1\leq s\leq \lfloor t/2\rfloor $. \\ \\
We then consider the case when $Z(G)\leq \Omega_{n-1}(G)$. If for some $1\leq i\leq 2s$, we have $[a_{i},b]=a^{p^{n-1}\alpha}$ then for a suitable choice of $j$ and $\beta$ we obtain $[a_{j}b^{\beta},b_{i}]=1$. Hence we may assume without loss of generality that we have chosen our generators such that $[a_{1},b]=\ldots =[a_{2s},b]$. By our assumption that $Z(G)\leq \Omega_{n-1}(G)$ we can now furthermore assume that $t>2s$ and that $[a_{2s+1},b]=b^{p^{n-1}}$. For this situation we thus get a group
$B(n,t,s)$ of order $p^{n+t}$ with presentation 
    $$\begin{array}{l}
                    a_{1}^{p}=\cdots =a_{t}^{p}=b^{p^{n}}=1, \\
                    \mbox{} [a_{2i-1},a_{2i}]=a^{p^{n-1}}\mbox{ for }i=1,\ldots ,s, \\
           \mbox{}        [a_{i},a_{j}]=1\mbox{ otherwise for }1\leq i<j\leq t, \\
      \mbox{}            [a_{1},b]=\ldots =[a_{2s},b]=[a_{2s+2},b]=\ldots =[a_{t},b]=1, \\
\mbox{}                 [a_{2s+1},b]=a^{p^{n-1}}.
\end{array}$$    
For such a group to be abelian we must have $n\geq 3$ as for $n=2$ we have $Z(G)^{p}=1$. Notice that for a fixed $n\geq 3$ and
$t\geq 1$ we get $\lfloor (t+1)/2\rfloor $ such groups, namely $B(n,t,s)$ for $0\leq s\leq \lfloor (t-1)/2\rfloor $. \\ \\
We are now ready for our classification. \\ \\
{\bf Groups of order $p^{4}$}. The only possible types for non-abelian groups of rank greater than $2$ is $(1,1,2)$. Apart from the
$5$ abelian groups we thus have 
                $$A(2,2,1)=\langle a,b,c:\,a^{p}=1,\, b^{p}=1,\, c^{p^{2}}=1,\,[a,b]=c^{p},\,[a,c]=[b,c]=1\rangle,$$ 
and
     $$G(3,1,2)=\langle a,b:\,a^{p^{3}}=b^{p}=1,\,[a,b]=a^{p^{2}}\rangle,$$
giving us in total $7$ groups. \\ \\
{\bf Groups of order $p^{5}$}. Here the possible types for non-abelain groups of rank greater than $2$ are $(1,1,1,2)$, $(1,1,3)$ and $(1,2,2)$.  The first two are covered by the discussion above, so we turn our attention to the last one. Let $G$ be a group of type $(1,2,2)$. As we mentioned in the introduction, we know from the structure theory of powerfully nilpotent groups that we can pick generators $a,b,c$ of order $p,p^{2},p^{2}$ for $G$ such that $G=\langle a\rangle\langle b\rangle\langle c\rangle$ and such that
we have a powerfully central series
     $$G>\langle b,c\rangle >\langle b^{p},c\rangle >\langle b^{p},c^{p}\rangle >\langle c^{p^{2}}\rangle >1.$$
Thus $[a,b]\in \langle b^{p},c\rangle^{p}=\langle c^{p}\rangle$ and $[a,c], [b,c]\in \langle b^{p},c^{p}\rangle^{p}=1$. Thus if the group is non-abelian we can
furthermore assume we have chosen the generators such that $[a,b]=c^{p}$. This gives us a unique powerfully nilpotent group with presentation
       $$E_{1}=\langle a,b,c:\,a^{p}=1,\, b^{p^{2}}=1,\, c^{p^{2}}=1,\,[a,b]=c^{p},\,[a,c]=[b,c]=1\rangle.$$ 
Apart from this there are $7$ abelian groups and we also have 
         $$G(3,2,2)=\langle a,b:\,a^{p^{3}}=b^{p^{2}}=1,\, [a,b]=a^{p^{2}}\rangle,$$
        $$G(4,1,3)=\langle a,b:\,a^{p^{4}}=b^{p}=1,\,[a,b]=a^{p^{3}}\rangle,$$
       $$A(2,3,1)=\langle a,b,c,d:\,a^{p}=b^{p}=c^{p}=d^{p^{2}}=1,\,[a,b]=d^{p},\,[a,c]=[a,d]=[b,c]=[b,d]=[c,d]=1\rangle,$$
      $$A(3,2,1)=\langle a,b,c:\,a^{p}=b^{p}=c^{p^{3}}=1,\, [a,b]=c^{p^{2}}\,[a,c]=[b,c]=1\rangle$$
and
      $$B(3,2,0)=\langle a,b,c:\, a^{p}=b^{p}=c^{p^{3}}=1,\,[a,c]=c^{p^{2}},[a,b]=[b,c]=1\rangle.$$
This gives us in total $13$ groups. \\ \\
{\bf Groups of order $p^{6}$.} Now the possible types for non-abelian groups of rank greater than $2$ are $(1,1,1,1,2)$, $(1,1,1,3)$, $(1,1,2,2)$, $(1,1,4)$,
$(1,2,3)$ and $(2,2,2)$. Those that are not covered by the general examples above are $(2,2,2)$, $(1,2,3)$ and $(1,1,2,2)$. Let $G$ be a group of type $(2,2,2)$. We know from the general structure theory that we can find generators $a,b,c$ of order $2$ such that $G=\langle a\rangle \langle b\rangle \langle c\rangle$ and where we have a powerfully central chain
      $$G>\langle a^{p},b,c\rangle>\langle a^{p},b^{p},c\rangle > \langle a^{p},b^{p},c^{p}\rangle >\langle b^{p},c^{p}\rangle >\langle c^{p}\rangle >1.$$
In particular $[a,c]=[b,c]\in G^{p^{2}}=1$ and $[a,b]\in \langle c^{p}\rangle$. There is thus a unique non-abelian group of type $(2,2,2)$, namely
       $$E_{2}=\langle a,b,c:\,a^{p^{2}}=b^{p^{2}}=c^{p^{2}}=1,\,[a,b]=c^{p},\,[a,c]=[b,c]=1\rangle.$$
We next turn to groups of type $(1,2,3)$. Here we can choose the generators $a,b,c$ of order $p,p^{2},p^{3}$ such that one of the following chains is
central
             $$G>\langle b,c\rangle >\langle b^{p},c\rangle > \langle b^{p},c^{p}\rangle > \langle c^{p}\rangle >\langle c^{p^{2}}\rangle >1$$
            $$G>\langle c,b\rangle >\langle c^{p},b\rangle > \langle c^{p},b^{p}\rangle > \langle c^{p^{2}},b^{p}\rangle >\langle c^{p^{2}}\rangle >1.$$
Let us first deal with former case. Notice that $[a,c],[b,c]\in G^{p^{2}}$. Also $[a,b]\in \langle b^{p},c\rangle^{p}=\langle c^{p}\rangle$. But as $[a,b]^{2}=[a^{p},b]=1$ we then must have that $[a,b]\in \langle c^{p^{2}}\rangle$ and thus $[G,G]\leq G^{p^{2}}=\langle c^{p^{2}}\rangle $ and $G$ is strongly powerful. It follows that $G^{p}\leq Z(G)$. This means that we get an alternating form $V\times V\rightarrow \mbox{GF}(p)$ that maps
$(xG^{p},yG^{p})\rightarrow \lambda$ where $[x,y]=c^{p^{2}\lambda}$. But $V$ has odd dimension and thus there is $x\in Z(G)\setminus G^{p}$. Thus $Z(G)\not =G^{p}$. We consider $3$ subcases. \\ \\
\underline{Case 1. $ Z(G)\leq \Omega_{1}(G)G^{p}=\langle a,b^{p},c^{p}\rangle$}. As $Z(G)\not\leq G^{p}$ we can pick our generators such that $a\in Z(G)$ and $G$ is direct product of $\langle a\rangle$ and $\langle b,c\rangle$. From our analysis of powerfully nilpotent groups of rank $2$ we can furthermore
pick $b,c$ such that we get a group with presentation
        $$E_{3}=\langle a,b,c\rangle:\,a^{p}=b^{p^{2}}=c^{p^{3}}=1,\,[b,c]=c^{p^{2}},\,[a,b]=[a,c]=1\rangle.$$
\underline{Case 2. $ Z(G)\not\leq \Omega_{1}(G)G^{p}=\langle a,b^{p},c^{p}\rangle$ but $Z(G)\leq \Omega_{2}(G)=\langle a,b,c^{p}\rangle$}. Here we
can assume that the generators have been picked such that $b\in Z(G)$. Here $G$ is the semidirect of $\langle b\rangle$ and $\langle a,c\rangle$. Our analysis of groups of rank $2$ shows that in this case we get a group with presentation 
      $$E_{4}=\langle a,b,c:\,a^{p}=b^{p^{2}}=c^{p^{3}}=1,\,[a,c]=c^{p^{2}},\,[a,b]=[b,c]=1\rangle.$$
\underline{Case 3. $ Z(G)\not\leq \Omega_{2}(G)=\langle a,b,c^{p}\rangle$}. In this case we can assume that our generators have been chosen such that 
$c\in Z(G)$. Without loss of generality we can also pick $a,b$ such that $[a,b]=c^{p^{2}}$. We now arrive at a group with presentation
    $$E_{5}=\langle a,b,c:\,a^{p}=b^{p^{2}}=c^{p^{3}}=1,\,[a,b]=c^{p^{2}},\,[a,c]=[b,c]=1\rangle.$$
Notice that if $G$ is strongly powerful then the series $G>G^{p}>G^{p^{2}}>1$ is powerfully central and we are in the situation that we have just been working with. We can thus assume that our group $G$ of type $(1,2,3)$ is not strongly powerful. We then have just seen that that the following chain is powerfully
central
     $$G>\langle c,b\rangle >\langle c^{p},b\rangle > \langle c^{p},b^{p}\rangle > \langle c^{p^{2}},b^{p}\rangle >\langle c^{p^{2}}\rangle >1.$$
Notice that in particular $[a,b],[b,c]\in G^{p^{2}}=\langle c^{p^{2}}\rangle$ and as $G$ is not strongly powerful we can't have $[a,c]\in G^{p^{2}}$. It follows
that $Z(G)\not\leq \Omega_{2}(G)=\langle a,b,c^{p}\rangle$, as otherwise there would be some $a^{\alpha}b^{\beta}c\in Z(G)$
giving the contradiction that $[a,c]\in G^{p^{2}}$. Now $[a,c]\in \langle c^{p},b\rangle^{p}$ and, as $[a,c]\not\in G^{p^{2}}$,  replacing $b$ by suitable $c^{pl}b$ we can assume that $[a,c]=b^{p}$. We now consider few subcases. \\ \\
\underline{Case 1. $ Z(G)\not\leq G^{p}$}. This means that there is an element of the form $a^{\alpha}b^{\beta}c^{\gamma}$ where at least one of  $\alpha$, $\beta$ and $\gamma$ is not divisible by $p$. If either $\alpha$ or $\gamma$ is not divisible by $p$ then after taking commtator with $a$ or $c$ we would get the contradiction that $[a,c]\in G^{p^{2}}$. We thus get a element $a^{\lambda}b^{\beta}c^{p\gamma}$ in $Z(G)$ and then also $\bar{b}=a^{\lambda}b^{\beta}\in Z(G)$ replacing $b$ by $\bar{b}$ we get a group with 
presentation 
     $$E_{6}=\langle a,b,c:\,a^{p}=b^{p^{2}}=c^{p^{3}}=1,\,[a,b]=[b,c]=1,\,[a,c]=b^{p}\rangle.$$
\underline{Case 2. $Z(G)=G^{p}$ and $\Omega_{2}(G)'=1$}. This means that $[a,b]=1$ and as $b\not\in Z(G)$ we then have
that $[b,c]=c^{p^{2}\lambda}$ for some $0<\lambda <p$. Replacing $b$ by $b^{\eta}$  and $a$ by $a^{\eta}$, where $\eta$ is the inverse of $\lambda$ modulo $p$, we can then assume that $[b,c]=c^{p^{2}}$ and $[a,c]=b^{p} $. We thus arrive at a group with presenation 
     $$E_{7}=\langle a,b,c:\,a^{p}=b^{p^{2}}=c^{p^{3}}=1\,[a,b]=1,\,[a,c]=b^{p},\,[b,c]=c^{p^{2}}.$$
\underline{Case 3. $Z(G)=G^{p}$ and $\Omega_{2}(G)'\not =1$}. This means that $[a,b]=c^{p^{2}\lambda}$ for some $\lambda$ coprime to $p$. Replacing $c$ by a suitable $ca^{s}$ we can assume that $[b,c]=1$. We thus arrive at a group with presentation
                           $$\langle a,b,c:\,a^{p}=b^{p^{2}}=c^{p^{3}}=1,\,[a,b]=c^{p^{2}\lambda},\,[a,c]=b^{p},\,[b,c]=1\rangle.$$
Now suppose the same group can be representated by $\bar{a},\bar{b},\bar{c}$ satisifying the same kind of relations. We want to
understand what the possible values for $\lambda$ in $\mbox{GF}(p)$ are.  In order to get same kind of relations we must to start with have
$$\bar{a}=a^{r}b^{ps}c^{p^{2}t},\ \bar{b}=a^{f}b^{g}c^{ph},\ \bar{c}=a^{i}b^{j}c^{k},$$
where $r,g$ and $k$ are coprime to $p$. For these to satisfy furthermore the three commutator relations, it is not difficult to
see that $p|f$ (as $[\bar{b},\bar{c}]=1$), $p|j$ (as $[\bar{a},\bar{c}]=\bar{b}^{p}$) and then $p|i$ (as $[\bar{b},\bar{c}]=1$).
We then get 
        $$[\bar{a},\bar{b}]=[a^{r},b^{g}]=c^{p^{2}rg\lambda},\, [\bar{a},\bar{c}]=[a^{r},c^{k}]=b^{prk}.$$
Thus for the relations to hold we need the following to hold in ${\mathbb K}=\mbox{GF}(p)$. 
        $$rk=g,\ rg\lambda=k\bar{\lambda}$$
where $\bar{\lambda}$ is the value replacing $\lambda$ for this new choice of basis. Thus 
       $$\bar{\lambda}=(rg/k)\lambda=r^{2}\lambda$$
We thus get that $\lambda$ can be replaced by any element in $({\mathbb K}^{*})^{2}\lambda$ and as $({\mathbb K}^{*})^{2}$ is of index $2$ we get exactly two groups with presentations
        $$E_{8}=\langle a,b,c:\,a^{p}=b^{p^{2}}=c^{p^{3}}=1,\, [a,b]=c^{p^{2}},\,[a,c]=b^{p},\,[b,c]=1\rangle$$
and 
       $$E_{9}=\langle a,b,c:\,a^{p}=b^{p^{2}}=c^{p^{3}}=1,\,[a,b]=c^{p^{2}\lambda},\,[a,c]=b^{p},\,[b,c]=1\rangle,$$
where $\lambda$ is any non-square in ${\mathbb K}^{*}$. \\ \\
We are now only left with groups of type $(1,1,2,2)$. In this case we pick our generators $a,b,c,d$ of orders $p,p,p^{2},p^{2}$
such that we have a central chain
   $$G>\langle b,c,d\rangle>\langle c,d\rangle>\langle c^{p},d\rangle>\langle c^{p},d^{p}\rangle>\langle d^{p}\rangle>1.$$
Notice that in particular we must have that any commutator with $d$ is in $\langle c^{p},d^{p}\rangle^{p}=1$ and thus
$d\in Z(G)$. Notice also that the group has real nilpotency class $2$ as $[G,G,G]\leq [G^{p},G]=[G,G]^{p}\leq G^{p^{2}}=1$. We consider again some cases. \\ \\
\underline{Case 1. $|Z(G)^{p}|=2$}. It is not difficult to see that one can assume that $c,d\in Z(G)$. For $G$ not to be abelian 
we then must have $[a,b]\not =1$. Without loss of generality we can then assume that $[a,b]=c^{p}$. We thus get the group
     $$E_{10}=\langle a,b,c,d:\,a^{p}=b^{p}=c^{p^{2}}=d^{p^{2}}=1,\,[a,b]=c^{p},\,[a,c]=[b,c]=[a,d]=[b,d]=[c,d]=1\rangle.$$ 
\underline{Case 2. $Z(G)^{p}=\langle d^{p}\rangle$ and $\Omega_{1}(G)$ is abelian}. Then $[a,b]=1$ and at least one of $[a,c],[b,c]$ is non-trivial. It is easy to see that without loss of generality we can assume that $[a,c]=1$ and $[b,c]=d^{p}$. We thus
get the group
    $$E_{11}=\langle a,b,c,d:\,a^{p}=b^{p}=c^{p^{2}}=d^{p^{2}}=1,\,[a,b]=1,\,[b,c]=d^{p},\,[a,c]=[a,d]=[b,d]=[c,d]=1\rangle.$$
\underline{Case 3. $Z(G)^{p}=\langle d^{p}\rangle$ and $\Omega_{2}(G)$ is not abelian.} Notice first that $G/Z(G)^{p}$ is
not abelian as otherwise we can choose the generators such that $[a,b]=[a,c]=d^{p}$ and $[b,c]=1$. But then $cb^{-1}\in Z(G)$ that gives the contradiction that $|Z(G)^{p}|=p$. It follows that we can pick our generators such that $[a,b]=c^{p}$,
$[a,c]=d^{p}$, and $[b,c]=1$. We thus get the group
     $$E_{12}=\langle a,b,c,d:\,a^{p}=b^{p}=c^{p^{2}}=d^{p^{2}}=1\,[a,b]=c^{p},\,[a,c]=d^{p},\,[a,d]=[b,c]=[b,d]=[c,d]=1\rangle.$$
Apart from the 11 groups $E_{2},\ldots ,E_{12}$ and the 11 abelian groups we have also the following groups of order $p^{6}$ (not writing down the trivial commutator relations)
\begin{eqnarray*}
  A(2,4,1) & = & \langle a,b,c,d,e:\,a^{p}=b^{p}=c^{p}=d^{p}=e^{p^{2}}=1,\,[a,b]=e^{p}\rangle, \\
   A(2,4,2) & = & \langle a,b,c,d,e:\,a^{p}=b^{p}=c^{p}=d^{p}=e^{p^{2}}=1,\,[a,b]=[c,d]=e^{p}\rangle, \\
 A(3,3,1) & = & \langle a,b,c,d:\,a^{p}=b^{p}=c^{p}=d^{p^{3}}=1,\,[a,b]=d^{p^{2}}\rangle, \\
 A(4,2,1) & = & \langle a,b,c:\,a^{p}=b^{p}=c^{p^{4}}=1,\,[a,b]=c^{p^{3}}\rangle, \\
 B(3,3,0) & = & \langle a,b,c,d:\,a^{p}=b^{p}=c^{p}=d^{p^{3}}=1,\,[a,d]=d^{p^{2}}\rangle, \\
B(3,3,1) & = & \langle a,b,c,d:\,a^{p}=b^{p}=c^{p}=d^{p^{3}}=1,\,[a,b]=d^{p^{2}},\,[a,d]=d^{p^{2}},\rangle, \\
B(4,2,0) & = & \langle a,b,c:\,a^{p}=b^{p}=c^{p^{4}}=1,\,[a,c]=c^{p^{3}}\rangle,
\end{eqnarray*}
and the following groups of rank $2$.
\begin{eqnarray*}
    G(3,3,2) & = & \langle a,b:\,a^{p^{3}}=b^{p^{3}}=1,\,[a,b]=a^{p^{2}}\rangle, \\
  G(5,1,4) & = & \langle a,b:\,a^{p^{5}}=b^{p}=1,\,[a,b]=a^{p^{4}}\rangle, \\
  G(4,2,2) & = & \langle a,b:\,a^{p^{4}}=b^{p^{2}}=1,\,[a,b]=a^{p^{2}}\rangle, \\
G(4,2,3) & = & \langle a,b:\,a^{p^{4}}=b^{p^{2}}=1,\,[a,b]=a^{p^{3}}\rangle.
\end{eqnarray*}
We have thus seen that there are 33 groups in total of order $p^{6}$.
\subsection{The case when $p=2$}
These are all the powerful groups of order less than $2^{6}$ that are strongly powerful. In particular all such groups
of exponent $4$ are abelian. The only groups that are of rank greater than $2$, exponent greater than $4$ and not covered by the groups of type $A(n,t,s)$ and $B(n,t,s)$ are groups of type $(1,2,3)$. \\ \\
{\bf Groups for order $16$}. The only non-abelian group is
      $$G(3,1,2)=\langle a,b:\,a^{8}=b^{2}=1,\,[a,b]=a^{4}\rangle.$$
Apart from this there are $5$ abelian groups and thus $6$ in total. \\ \\
{\bf Groups of order $32$}. Here the only non-abelian groups are
  \begin{eqnarray*}
                 G(3,2,2) & = & \langle a,b:\,a^{8}=b^{4}=1,\,[a,b]=a^{4}\rangle, \\
              G(4,1,3) & = & \langle a,b:\,a^{16}=b^{2}=1,\,[a,b]=a^{8}\rangle, \\
             A(3,2,1) & = & \langle a^{2}=b^{2}=c^{8}=1,\,[a,b]=c^{4}\rangle, \\
           B(3,2,0) & = & \langle a^{2}=b^{2}=c^{8}=1,\,[a,c]=c^{4}\rangle.
\end{eqnarray*}
Apart from this there are $7$ abelian groups giving a total of $11$ groups. \\ \\
{\bf Groups of order $64$}. The only exceptional groups that we need to consider are those of type $(1,2,3)$. Here we can pick
generators $a,b,c$ of orders $1,2,3$ such that $G=\langle a\rangle\langle b\rangle\langle c\rangle$ and we have a central
chain
      $$G>G^{p}>1.$$
 The only non-abelian groups here are $E_{3}$, $E_{4}$ ,$E_{5}$, $A(3,3,1)$, $A(4,2,1)$, $B(3,3,0)$, $B(3,3,1)$, $B(4,2,0)$ , $G(3,3,2)$, $G(5,1,4)$,
$G(4,2,2)$  and $G(4,2,3)$. These together with 11 abelian groups give in total 23 groups.

\end{document}